\documentclass[11pt,reqno]{amsart}
\usepackage{amsmath,amssymb,geometry,color,tikz-cd}
\usepackage[backref,pagebackref,pdftex,hyperindex]{hyperref}
\usepackage[alphabetic]{amsrefs}
\usepackage{amssymb}% see geometry.pdf on how to lay out the page. There's lots.
\usepackage{bbm}
\usepackage{enumitem}
\usepackage{upgreek}

\newtheorem{proposition}{Proposition}[section]
\newtheorem{theorem}[proposition]{Theorem}
\newtheorem{lemma}[proposition]{Lemma}
\newtheorem{corollary}[proposition]{Corollary}

\theoremstyle{definition}
\newtheorem{remark}[proposition]{Remark}
\newtheorem{definition}[proposition]{Definition}

\newtheorem{question}[proposition]{Question}

\title{Effective semi-ampleness of Hodge line bundles on curves}
\author{Chuyu Zhou}
\address{\'Ecole Polytechnique F\'ed\'erale de Lausanne (EPFL), MA C3 615, Station 8, 1015 Lausanne, Switzerland}
\email{chuyu.zhou@epfl.ch}
\date{} % delete this line to display the current date

\thanks{
		\newline
		\indent
		\emph{Competing interests}: The author declares none.
	}

\newcommand{\ord}{{\rm {ord}}}

\newcommand{\vol}{{\rm {vol}}}
\newcommand{\lct}{{\rm {lct}}}

\newcommand{\CM}{{\rm {CM}}}
\newcommand{\Hodge}{{\rm {Hodge}}}
\newcommand{\Mob}{{\rm {Mob}}}
\newcommand{\Fix}{{\rm {Fix}}}

\newcommand{\CY}{{\rm {CY}}}
\newcommand{\Bs}{{\rm {Bs}}}

\newcommand{\bC}{\mathbb{C}}

\newcommand{\bE}{\mathbb{E}}

\newcommand{\bN}{\mathbb{N}}

\newcommand{\bQ}{\mathbb{Q}}

\newcommand{\mD}{\mathcal{D}}

\newcommand{\mF}{\mathcal{F}}

\newcommand{\mL}{\mathcal{L}}
\newcommand{\mM}{\mathcal{M}}

\newcommand{\mO}{\mathcal{O}}

\newcommand{\mX}{\mathcal{X}}

\newcommand{\fB}{\mathbf{B}}

\begin{document}

\begin{abstract}
In this note, we prove effective semi-ampleness conjecture due to Prokhorov and Shokurov for a special case, more concretely, for $\bQ$-Gorenstein klt-trivial fibrations over smooth projective curves whose fibers are all klt log Calabi-Yau pairs of Fano type.
\end{abstract}

\maketitle
\tableofcontents

Throughout, we work over $\bC$.

\section{Introduction}\label{sec: 1}

Let $f: (X, \Delta)\to Z$ be a klt (resp. lc)-trivial fibration over a projective normal variety $Z$, i.e. $f$ is a projective and surjective morphism  with connected fibers and $(X,\Delta)$ is a klt (resp. lc) log pair such that $K_X+\Delta\sim_\bQ f^*D$ for some $\bQ$-divisor $D$ on $Z$. It is well known that  there exist two $\bQ$-divisors on $Z$, which are called discriminant part denoted by $B_Z$ and moduli part denoted by $M_Z$ such that $D\sim_\bQ K_Z+B_Z+M_Z$. As the names suggest, the divisor $B_Z$ measures the singularities of fibers while the divisor $M_Z$ measures the variation of fibers. Suppose $Z$ is an Ambro model (see \cite{Ambro04, Ambro05}), the effective semi-ampleness conjecture proposed by Prokhorov and Shokurov (\cite[Conjecture 7.13]{PS09}) says that, there should be a positive integer $m(d,r)$ depending only on the dimension $d$ and the torsion index $r$ of general fibers such that $mM_Z$ is base point free. We note here that the semi-ampleness of $M_Z$ is well known if $Z$ is a curve (e.g.  \cite{Ambro05, Floris14}).
The effective semi-ampleness is harder and only a few cases are known.
For example, O.Fujino \cite[Theorem 1.2]{Fujino03} shows that if general fibers are K3 surfaces or abelian varieties of dimension $d$, then one can take $m=19k$ and $m=k(d+1)$ respectively to achieve base point freeness, where $k$ is a weight associated to the Baily-Borel-Satake
compactification of the period domain. It is worth mentioning that the proof heavily depends on the existence of moduli space of fibers. The idea is that $M_Z$ appears due to the variation of fibers, it is apparent that one will have a good understanding of it once the moduli space of fibers exists.

Recently, a projective separated good moduli space has been constructed for Fano varieties with fixed invariants and K-stability due to works \cite{Jiang20, BLX19, Xu20, BX19, ABHLX20, BHLLX20, XZ20b, LXZ21}. It is the so called K-moduli space. Thus one can apply the idea to the klt-trivial fibrations whose fibers are parametrized by a K-moduli space. For a klt (resp. lc) trivial fibration $(X,\Delta)\to C$ where $C$ is a curve, we say that it is $\bQ$-Gorenstein if $K_{X/C}$ is $\bQ$-Cartier. The following theorem is the main result of this paper.

\begin{theorem}\label{main}
Fix positive integers $d$ and $r>1$. Let $f: (X,\Delta)\to C$ be a 
$\bQ$-Gorenstein klt-trivial fibration over a smooth projective curve $C$ such that
\begin{enumerate}
\item for any closed point $t\in C$, the fiber $X_t$ is a $\bQ$-Fano variety of dimension $d$, i.e. $X_t$ admits klt singularities and $-K_{X_t}$ is ample,
\item $r\Delta_t$ is a Weil divisor on $X_t$ and no component of $\Delta$ is contained in a fiber,
\item for any closed point $t\in C$, the fiber $(X_t, \Delta_t)$ is klt with $K_{X_t}+\Delta_t\sim_\bQ 0$,
\end{enumerate}
Write $K_X+\Delta=f^*(K_C+\Lambda_C)$.
Then there exists a positive integer $m(d, r)$ depending only on $d$ and $r$ such that $m\Lambda_C$ is base point free.
\end{theorem}

Here $\Lambda_C$ is the moduli part of the klt-trivial fibration, also called Hodge line bundle (which is in fact a $\bQ$-line bundle) in this case. We roughly talk here about the idea of the proof which will be treated more carefully later. As $(X_t, \Delta_t)$ is a klt log Calabi-Yau with $-K_{X_t}$ being ample, the pair $(X_t, (1-\epsilon)\Delta_t)$ is a K-stable log Fano pair for any $0<\epsilon\ll 1$. Via K-stability, one can construct a K-moduli space $M$ parametrizing the fibers with a Hodge line bundle $\Lambda_{\Hodge}$ on $M$.  For any klt-trivial fibration as in Theorem \ref{main}, there is a morphism $C\to M$ such that $\Lambda_C$ is the pullback of the Hodge line bundle on $M$. Thus one can read the information of $\Lambda_C$ from $\Lambda_{\Hodge}$. Note that in Theorem \ref{main}, we require that every fiber is a \textit{klt} log Calabi-Yau pair, which is a key point in the proof. In the last section, we will discuss about a more general setting.

\noindent
\subsection*{Acknowledgement}
   The author would like to thank Chen Jiang, Junpeng Jiao, Yuchen Liu, Zhan Li,  Zsolt Patakfalvi, Javier Carvajal Rojas and Roberto Svaldi for helpful discussions. The author is grateful to Ziquan Zhuang for answering his questions on the paper \cite{XZ20b}. The author is supported by grant European Research Council (ERC-804334).

\section{Preliminaries}\label{sec: 2}

In this section,  we give a quick introduction of the concepts of  K-stability, CM-line bundle, and Hodge line bundle. We say $(X,\Delta)$ is a log pair if $X$ is a projective normal variety and $\Delta$ is an effective $\bQ$-divisor on $X$ such that $K_X+\Delta$ is $\bQ$-Cartier.  The log pair $(X,\Delta)$ is called log Fano if it admits klt singularities and $-(K_X+\Delta)$ is ample; if $\Delta=0$, we just say $X$ is a $\bQ$-Fano variety. For the concepts of klt and lc singularities, please refer to \cite{KM98,Kollar13}.

\subsection{K-stability}

Let $(X,\Delta)$ be a log pair. Suppose $f: Y\to X$ is a proper birational mophism between normal varieties and $E$ is a prime divisor on $Y$, we say that $E$ is a prime divisor over $X$ and we define the following invariant
$$A_{X,\Delta}(E):=\ord_E(K_Y-f^*(K_X+\Delta)), $$
which is called the log discrepancy of $E$ associated to the pair $(X,\Delta)$.
If $(X,\Delta)$ is a log Fano pair, we define the following invariant
$$S_{X,\Delta}(E):=\frac{1}{\vol(-K_X-\Delta)}\int_0^\infty \vol(-f^*(K_X+\Delta)-xE){\rm{d}}x .$$
Denote $\beta_{X,\Delta}(E):=A_{X,\Delta}(E)-S_{X,\Delta}(E)$. By the works \cite{Fuj19, Li17, BJ20}, one can define K-stability of a log Fano pair by the following way.

\begin{definition}
Let $(X,\Delta)$ be a log Fano pair. 
\begin{enumerate}
\item We say that $(X,\Delta)$ is K-semistable if $\beta_{X,\Delta}(E)\geq 0$ for any prime divisor $E$ over $X$.
\item We say that $(X,\Delta)$ is K-polystable if it is K-semistable and any prime divisor $E$ over $X$ satisfying $\beta_{X,\Delta}(E)= 0$ is of product type (see \cite[Section 3.2]{Fuj17}).
\item We say that $(X,\Delta)$ is uniformly K-stable if there is a positive rational number $0< \epsilon\ll 1$ such that
$\beta_{X,\Delta}(E)\geq \epsilon S_{X,\Delta}(E)$ for any prime divisor $E$ over $X$.
\end{enumerate}

\end{definition}

\subsection{CM-line bundle and Hodge line bundle}

In this section, we fix an integer $r>1$ and a family $\pi: (\mX, \frac{1}{r}\mD)\to B$ of relative dimension $d$ satisfying the following conditions:
\begin{enumerate}
\item $\pi$ and $\pi|_\mD$ are flat projective morphisms, and $\mD$ is a Weil divisor on $\mX$,
\item $-K_{\mX/B}$ is a relatively ample $\bQ$-Cartier divisor on $\mX$ and  $\mD\sim_{\bQ,B} -rK_{\mX/B}$,\item $(\mX, \frac{1}{r}\mD+\mX_t)$ is slc for every closed point $t\in B$.
\end{enumerate}

Let $\mL$ be a relatively ample $\bQ$-line bundle on $\mX$. We have the following Knudsen-Mumford expansion (see \cite{KM76}):
$$\det(\pi_*(k\mL))=\lambda_{d+1}^{\binom{k}{d+1}}\otimes\lambda_d^{\binom{k}{d}} \otimes...\otimes\lambda_1^{\binom{k}{1}}\otimes\lambda_0,$$
$$\det({\pi|_\mD}_*(k\mL|_\mD))=\tilde{\lambda}_d^{\binom{k}{d}}\otimes...\otimes\tilde{\lambda}_0, $$
where $\lambda_i, i=0,1,...,d+1,$ and $\tilde{\lambda}_j, j=0,1,...,d,$ are $\bQ$-line bundles on $B$. 
By the flatness of $\pi$ and $\pi|_\mD$, we know that $h^0(\mX_t, k\mL_t)$ and $h^0(\mD_t, k{\mL_t}|_{\mD_t})$ are independent of $t$. We write
$$h^0(\mX_t, k\mL_t)=a_0k^d+a_1k^{d-1}+o(k^{d-1})\quad and \quad  h^0(\mD_t, k{\mL_t}|_{\mD_t})=\tilde{a}_0k^{d-1}+o(k^{d-1}).$$
It is not hard to see
$$a_0=\frac{\mL_t^d}{d!}, 
a_1=\frac{-K_{\mX_t}\mL_t^{d-1}}{2(d-1)!}, 
\tilde{a}_0=\frac{\mL^{d-1}_t\mD_t}{(d-1)!}.$$

\begin{definition}\label{def: cm}
We define a $\bQ$-line bundle $\lambda_{\CM,c}$ (i.e. CM-line bundle) on $B$ with respect to the polarization $\mL$ as follows:
$$\mM_1:=\lambda_{d+1}^{\frac{2a_1}{a_0}+d(d+1)}\otimes\lambda_d^{-2(d+1)}\quad and \quad \mM_2:=\lambda_{d+1}^{\frac{\tilde{a}_0}{a_0}} \otimes \tilde{\lambda}_d^{-(d+1)},$$
$$\lambda_{\CM,c}:= \mM_1-\frac{c}{r}\mM_2.$$
\end{definition}

\begin{definition}\label{def: hodge}
We define a $\bQ$-line bundle $\Lambda_B$ (i.e. Hodge line bundle) on $B$ as follows:
$$K_{\mX/B}+\frac{1}{r}\mD= \pi^*\Lambda_B. $$
\end{definition}

The CM-line bundle satisfies the base change property, see \cite[Lemma 3.5, Proposition 3.8]{CP21}. When $B$ is a projective normal variety,
by Riemann-Roch formula, cf \cite[Appendix]{CP21}, we have
\begin{enumerate}
\item $\lambda_{d+1}=\pi_*(\mL^{d+1}), $
\item $\lambda_d=\frac{d}{2}\pi_*(\mL^{d+1})+\frac{1}{2}\pi_*(-K_{\mX/B}\mL^d),$
\item $\tilde{\lambda}_d=\pi_*(\mL^d\mD). $
\end{enumerate}

\begin{remark}\label{cases on cm}
Notation as above, if $B$ is a projective normal base, we have the following simple computations:
\begin{enumerate}
\item $\mM_1= (d+1)\pi_*(\mL^dK_{\mX/B})-\frac{d(K_{\mX_t}\mL_t^{d-1})}{\mL_t^d}\pi_*\mL^{d+1}.$
\item $\mM_2=-(d+1)\pi_*(\mL^d\mD)-\frac{d(\mD_t\mL_t^{d-1})}{\mL_t^d}\pi_*\mL^{d+1}.$
\item $\lambda_{\CM,c}:=(d+1)\pi_*(\mL^d(K_{\mX/B}+\frac{c}{r}\mD))-\frac{d(K_{\mX_t}+\frac{c}{r}\mD_t)\mL_t^{d-1}}{\mL_t^d}\pi_*\mL^{d+1}.$
\item If $c\in [0,1)$ and $\mL=-(K_{\mX/B}+\frac{c}{r}\mD)$, then
$\lambda_{\CM,c}= -\pi_*\mL^{d+1}.$
\item If $c>1$ and $\mL=K_{\mX/B}+\frac{c}{r}\mD$, then
$\lambda_{\CM,c}= \pi_*\mL^{d+1}.$
\item If $c=1$, then
$\lambda_{\CM,1}= (d+1)\pi_*(\mL^d(K_{\mX/B}+\frac{1}{r}\mD)).$
\end{enumerate}
\end{remark}

\begin{lemma}\label{cm vs hodge}
We have the following relation between Hodge line bundle and CM-line bundle:
$$\lambda_{\CM,1}=(d+1)\mL_t^d\cdot\Lambda_B. $$
\end{lemma}
\begin{proof}
The proof is a combination of Definition \ref{def: hodge} and Remark \ref{cases on cm} (6).
\end{proof}

\section{K-moduli space with a Hodge line bundle}\label{sec: 3}

In this section, we fix two integers $d$ and $r$ and consider a set $\mF$ of log pairs satisfying that $(X, \frac{1}{r}D)\in \mF$ if and only if the following conditions are satisfied:
\begin{enumerate}
\item $X$ is a $\bQ$-Fano variety of dimension $d$ and $\vol(-K_X)=v$,
\item $D\sim_\bQ -rK_X$ is a Weil divisor on $X$,
\item $(X,\frac{1-\epsilon}{r}D)$ is K-semistable for any rational $0<\epsilon\ll 1$.
\end{enumerate}
By \cite[Theorem 5.4]{Zhou21b} and \cite[Lemma 5.3]{Zhou21b}, we know that 
\begin{enumerate}
\item the set $\mF$ lies in a log bounded family,
\item there exists a positive rational number $0<\epsilon_0<1$ depending only on $d$ and $r$ such that if $(X,\frac{1}{r}D)\in \mF$, the pair $(X, \frac{c}{r}D)$ is K-semistable for all rational $c\in [1-\epsilon_0,1)$.
\end{enumerate}
Note that every log pair $(X,\frac{1}{r}D)\in \mF$ admits log canonical singularities. Fix a rational number $c_0\in [1-\epsilon_0, 1)$, by \cite{BLX19,Xu20}, there exists an Artin stack of finite type, denoted by $\mM^K_{d,v,r}$, parametrizing 
log pairs $(X,\frac{c_0}{r}D)$ with $(X,\frac{1}{r}D)\in \mF$, where $\mM^K_{d,v,r}(S)$ consists of all families $(\mX, \frac{c_0}{r}\mD)\to S$ satisfying the following conditions:
\begin{enumerate}
\item $\mX\to S$ is proper and flat,
\item $\mD$ is a K-flat family of divisors on $\mX$ (see \cite{Kollar19}),
\item $-K_{\mX/S}-\frac{c_0}{r}\mD$ is $\bQ$-Cartier,
\item all geometric fibers are contained in $\mF$ after one changes the coefficient $\frac{c_0}{r}$ to be $\frac{1}{r}$.
\end{enumerate}
Note that the Artin stack $\mM^K_{d,v,r}$ does not depend on the choice of $c_0\in [1-\epsilon_0, 1)$ as K-stability does not change as $c_0$ varies in $[1-\epsilon_0, 1)$.
By \cite{BX19, ABHLX20, XZ20b, LXZ21},  $\mM^K_{d,v,r}$ admits a projective separated good moduli space, denoted by $M^K_{d,v,r}$. Note that by our notation, a closed point in $M^K_{d,v,r}$ corresponds to a log pair $(X,\frac{1}{r}D)\in \mF$ such that $(X,\frac{1-\epsilon}{r}D)$ is K-polystable for rational $0<\epsilon\ll 1$.  For convenience, we still use $\mF$ to denote the scheme which parametrizes the log pairs in $\mF$. Consider the universal family over $\mF$, denoted by $(\mX, \frac{1}{r}\mD)\to \mF$, with the polarization $\mL:=-K_{\mX/\mF}$.  By Definition \ref{def: cm}, one can define the CM-line bundle  on the base $\mF$ as follows
$$\lambda_{\CM,c}=\mM_1-\frac{c}{r}\mM_2 ,$$ 
where $c\in [1-\epsilon_0, 1+\epsilon_0]$, and $\mM_1$ and $\mM_2$ are $\bQ$-line bundles on $\mF$ defined by Knudsen-Mumford expansion. 
By \cite{LWX18b, CP21}, $\lambda_{\CM,c}$ descends to a $\bQ$-line bundle on $M^K_{d,v,r}$ for each rational $c\in [1-\epsilon_0, 1)$, denoted by $\Lambda_{\CM,c}$. 

\begin{lemma}\label{pullback lemma}
For any rational $c\in [1, 1+\epsilon_0)$, the $\bQ$-line bundle $\lambda_{\CM,c}$ can be descended to a $\bQ$-line bundle on $M^K_{d,v,r}$, also denoted by $\Lambda_{\CM,c}$.
In particular, the $\bQ$-line bundle $\frac{1}{(d+1)v}\lambda_{\CM,1}$ can be descended to a $\bQ$-line bundle on $M^K_{d,v,r}$, denoted by $\Lambda_{\Hodge}$. Moreover, for any morphism from a projective normal variety $B\to \mF$, let $\pi_B: (\mX_B,\frac{1}{r}\mD_B)\to B$ be the induced family via pulling back the universal family over $\mF$ and $\Lambda_B$ the Hodge line bundle on $B$, i.e. $K_{\mX_B/B}+\mD_B=\pi_B^*\Lambda_B$, then we have $\Lambda_B=f^*\Lambda_{\Hodge}$, where $f: B\to M^K_{d,v,r}$ is the induced morphism. 
\end{lemma}

\begin{proof}
To show that $\lambda_{\CM,c}$ for a rational $c\in [1, 1+\epsilon_0]$ can be descended to $M^K_{d,v,r}$, it suffices to show that both $\mM_1$ and $\mM_2$ can be descended to $M^K_{d,v,r}$. Take two different rational numbers $c', c''\in [1-\epsilon_0,1)$. As we have mentioned, both $\lambda_{\CM,c'}$ and $\lambda_{\CM,c''}$ can be descended to $M^K_{d,v,r}$, thus $\lambda_{\CM,c'}-\lambda_{\CM,c''}$, which is proportional to $\mM_2$, can be descended to $M^K_{d,v,r}$. Therefore, $\mM_1=\lambda_{\CM,c}+\frac{c}{r}\mM_2$ can also be descended to $M^K_{d,v,r}$. The last statement follows from Lemma \ref{cm vs hodge}.
\end{proof}

\section{Positivity of the Hodge line bundle}\label{sec: 4}

In the previous section we see that there is a projective separated good moduli space parametrizing all pairs $(X,\frac{1}{r}D)\in \mF$ such that $(X,\frac{1-\epsilon}{r}D)$ is K-polystable for sufficiently small rational $0<\epsilon\ll 1$. Moreover, there is a $\bQ$-line bundle $\Lambda_{\Hodge}$ (called Hodge line bundle) satisfying the pullback property as in Lemma \ref{pullback lemma}. In this section we study the positivity of the Hodge line bundle $\Lambda_{\Hodge}$.

\begin{theorem}\label{positivity}
Let $M$ be a proper subvareity of $M^K_{d,v,r}$ such that $M$ contains a closed point which corresponds to a klt log Calabi-Yau pair $(X,\frac{1}{r}D)\in \mF$, then the restriction of the Hodge line bundle $\Lambda_{\Hodge}$ on $M$ is big and nef.
\end{theorem}

\begin{proof}
We first show that $\Lambda_{\Hodge}$ itself is nef.  Take a sequence of strictly increasing rational numbers $\{c_i\}_{i\in \bN}$ such that $c_i\in [1-\epsilon_0, 1)$ and $\lim_{i\to \infty}c_i=1$. By the notation of Section \ref{sec: 3}, $\{\Lambda_{\CM,c_i}=M_1-\frac{c_i}{r}M_2\}_{i\in \bN}$ is a sequence of CM-line bundles on $M^K_{d,v,r}$, thus $\Lambda_{\CM,c_i}$ is ample for each $i$ by \cite{XZ20b}. By the proof of Lemma \ref{pullback lemma}, $\Lambda_{\Hodge}$ is equal to $\lim_{i\to \infty} \Lambda_{\CM,c_i}$ up to a positive multiple, thus the nefness is concluded. 

We turn to the bigness. By Lemma \ref{cm vs hodge} and recall the notation in Lemma \ref{pullback lemma} we see
$$\Lambda_{\CM,1}= (d+1)v\Lambda_{\Hodge}=\frac{1}{2}(\Lambda_{\CM,1-\epsilon}+\Lambda_{\CM,1+\epsilon})$$ 
for any rational $0<\epsilon\ll 1$.
By \cite{XZ20b} we know that $\Lambda_{\CM,1-\epsilon}$ is ample. It suffices to show that the restriction of $\Lambda_{\CM,1+\epsilon}$ on $M$ is pseudo-effective. Let $(\mX', \frac{1}{r}\mD')\to C$ be a  family over a smooth projective curve $C$ such that  
\begin{enumerate}
\item $(\mX', \frac{1-c_0}{r}\mD')\to C$ is contained in $\mM^K_{d,v,r}(C)$ (see Section \ref{sec: 3}),
\item for general $t\in C$, the fiber $(\mX'_t, \frac{1}{r}\mD'_t)$ is klt,
\item the image of the induced morphism $C\to M^K_{d,v,r}$ is contained in $M$.
\end{enumerate}
We denote $f: C\to M$ the induced morphism in the moduli sense, and consider the family
$(\mX', \frac{1+\epsilon}{r}\mD')\to C$ with the polarization $\mL_C:=-K_{\mX'/C}$.
 By ACC of log canonical thresholds (\cite{HMX14}), the general fibers of the family $(\mX', \frac{1+\epsilon}{r}\mD')\to C$ are still log canonical.
By \cite[Theorem 1.10]{Fujino18}, we know that $\pi_*\mO_{\mX'}(m(K_{\mX'/C}+\frac{1+\epsilon}{r}\mD'))$ is nef, which implies that $K_{\mX'/C}+\frac{1+\epsilon}{r}\mD'$ is also nef. By Definition \ref{def: cm}, we have
$$\deg(f^*\Lambda_{\CM,1+\epsilon})=(d+1)\mL_C^d(K_{\mX'/C}+\frac{1+\epsilon}{r}\mD')-\frac{d(K_{\mX'_t}+\frac{1+\epsilon}{r}\mD'_t)\mL_t^{d-1}}{\mL_t^d}\mL_C^{d+1}. $$ 
In our case, $\mL_C=\frac{1}{\epsilon}(K_{\mX'/C}+\frac{1+\epsilon}{r}\mD')+\pi^*H$ for some divisor $H$ on $C$. A direct computation implies that
$$\deg (f^*\Lambda_{\CM, 1+\epsilon})=\frac{1}{\epsilon^d}(K_{\mX'/C}+\frac{1+\epsilon}{r}\mD')^{d+1}\geq 0. $$
Therefore, the restriction of $\Lambda_{CM,1+\epsilon}$ on $M$ is pseudo-effective by \cite{BDPP13}. The proof is finished.
\end{proof}

\section{Proof of the main result}

In this section, we prove the main result in Section \ref{sec: 1}. We start by the following well-known lemma.

\begin{lemma}\label{kps}
Let d be a positive integer number, I a finite positive rational set contained in $[0,1]$. Consider the set $\{(X,D) \}$ satisfying the following conditions
\begin{enumerate}
\item $X$ is a $\bQ$-Fano variety of dim d,
\item the coefficients of D are contained in I,
\item $(X,D)$ is a klt Calabi-Yau pair.
\end{enumerate}
Then there is a positive number $\eta(d,I)$ such that $(X,(1-a)D)$ is K-polystable (even uniformly K-stable) for any rational $0<a <\eta$.
\end{lemma}

\begin{proof}
By ACC of log canonical thresholds, we know that there is a positive rational number $\epsilon$ such that $(X,D)$ is $\epsilon$-lc for any element $(X,D)$ in the set (e.g. \cite[Lemma 2.48]{Birkar19}), thus the set lies in a log bounded family due to \cite{HMX14} or \cite{Birkar19,Birkar21}. So for each $(X,D)$ in the set, one can find an uniform positive integer number $m$ and a very ample line bundle $A\sim -mK_X$ on $X$ satisfying conditions in \cite[Theorem 1.8]{Birkar21}. By the theorem, there is a positive real number $\delta(d,I)$ which only depends on $d$ and $I$ such that $\lct(X,D;|-K_X|_\bQ)\geq \delta$. That means, for each prime divisor $E$ over $X$, we have
$$A_{X,D}(E)-\delta T_X(E)\geq 0, $$
where $T_X(E)$ is the pseudo-effective threshold of $E$ with respect to $-K_X$.
A simple computation implies $S_{X, (1-a)D}(E)=aS_X(E)$.
Take $\eta:=\frac{\delta(d,I)}{2}$, then for any rational $0<a<\eta$,
$$\beta_{X,(1-a)D}(E)=A_{X,(1-a)D}(E)-a S_X(E)\geq A_{X,D}(E)-a T(E)\geq \frac{\delta T(E)}{2}\geq \frac{\delta}{2a}S_{X, (1-a)D}(E), $$
which deduces that $(X,(1-a)D)$ is uniformly K-stable for any rational $0<a<\eta$.
\end{proof}

By the above lemma we see that, for any fiber $(X_t,\Delta_t)$ of the family as in Theorem \ref{main}, the log pair $(X_t, (1-\epsilon)\Delta_t)$ is K-polystable for any rational $0<\epsilon\ll 1$. 

\begin{lemma}
All fibers of all families as in Theorem \ref{main} lie in a log bounded family.
\end{lemma}

\begin{proof}
This is deduced by applying \cite[Corollary 1.7]{HMX14}.
\end{proof}

By the above boundedness, we see that there are only finitely many volumes appearing for all the fibers $X_t$ of all the families as in Theorem \ref{main}. From now on, we can assume all the fibers admit volume $v$, where $v$ is a fixed positive rational number. Recall that in Section \ref{sec: 3} and Section \ref{sec: 4}, we have constructed a projective separated good moduli space $M^K_{d,v,r}$ with the Hodge line bundle $\Lambda_{\Hodge}$ admitting certain positvity, such that any fiber of any family as in Theorem \ref{main} corresponds to a closed point in $M^K_{d,v,r}$ (by Lemma \ref{kps}). We are ready to prove Theorem \ref{main}.

\begin{proof}[Proof of Theorem \ref{main}]
For the given family $(X,\Delta)\to C$, there is an induced morphism $f: C\to M^K_{d,v,r}$ such that $\Lambda_C=f^*\Lambda_{\Hodge}$. As all fibers are klt log Calabi-Yau pairs, by Theorem \ref{positivity} we may assume $M^{K}_{d,v,r}$ is irreducible and $\Lambda_{\Hodge}$ is big and nef. 
We first show that $f(C)\cap \fB(\Lambda_{\Hodge})=\emptyset$, where $\fB(\Lambda_{\Hodge})$ is the stable base locus of $\Lambda_{\Hodge}$. By \cite[Theorem 1.4]{Birkar17}, we have the following characterization for the augmented sable base locus $\fB_+(\Lambda_{\Hodge})$:
$$\fB_+(\Lambda_{\Hodge})=\bE(\Lambda_{\Hodge}), $$
where 
$$\bE(\Lambda_{\Hodge})=\cup_W \{\textit{$W\subset X$ is a proper subvariety such that ${\Lambda_{\Hodge}}|_W$ is not big}\}. $$
Let $z\in f(C)$ be a closed point and $W\subset M^K_{d,v,r}$ any proper subvariety containing $z$, again by Theorem \ref{positivity}, we see that $\Lambda_{\Hodge}|_W$ is big and nef. Therefore $z$ is not contained in $\bE(\Lambda_{\Hodge})$, which implies that $f(C)\cap \fB_+(\Lambda_{\Hodge})=\emptyset$. As $\fB(\Lambda_{\Hodge})\subset \fB_+(\Lambda_{\Hodge})$,  we also have $f(C)\cap \fB(\Lambda_{\Hodge})=\emptyset$.

Choose a sufficiently divisible $N$ such that $\fB(\Lambda_{\Hodge})=\Bs|N\Lambda_{\Hodge}|$.
Then take a partial resolution of $\Bs|N\Lambda_{\Hodge}|$, denoted by  $r: \tilde{M}\to M^{K}_{d,v,r}$, such that 
$$r^* N\Lambda_{\Hodge}\sim A+F,$$ 
where $A$ is a free line bundle and $\tilde{M}\setminus F\cong M^K_{d,v,r}\setminus \fB(\Lambda_{\Hodge})$. As $f(C)\cap \fB(\Lambda_{\Hodge})=\emptyset$,
the morphism $f: C\to M^K_{d,v,r}$ lifts to a morphism $\tilde{f}: C\to \tilde{M}$ such that $\tilde{f}(C)\cap F=\emptyset$. Thus $N\Lambda_C\sim\tilde{f}^* A$ and one can choose a positive integer $m$ depending only on $(M^K_{d,v,r}, \Lambda_{\Hodge})$ (therefore depending only on $d$ and $r$) such that $m\Lambda_C$ is base point free. The proof is finished.
\end{proof}

\begin{remark}
The same proof applies to the case where the base of the family in Theorem \ref{main} has a higher dimension.
\end{remark}

\section{Further exploration: Semi-ampleness of Hodge line bundles}

\subsection{A more general setting}

It is easy to ask whether Theorem \ref{main} still holds if one allows  strictly lc fibers for the family, as this is also a special case of the effective semi-ampleness conjecture of Prokhorov-Shokurov. We put the following question which is a little more general than Theorem \ref{main}. 

\begin{question}\label{lc fiber}
Fix positive integers $d$ and $r>1$. Let $f: (X,\Delta)\to C$ be a 
$\bQ$-Gorenstein lc-trivial fibration over a smooth projective curve $C$ such that

\begin{enumerate}
\item for any closed point $t\in C$, the fiber $X_t$ is a $\bQ$-Fano variety of dimension $d$,
\item $r\Delta_t$ is a Weil divisor on $X_t$ and no component of $\Delta$ is contained in a fiber,
%\item for general closed point $t\in C$, the log pair $(X_t, \frac{1}{r}\Delta_t)$ is klt log Calabi-Yau,
\item for every closed point $t\in C$, $K_{X_t}+\Delta_t\sim_\bQ 0$ and the log pair $(X_t, (1-\epsilon)\Delta_t)$ is K-semistable for any rational $0< \epsilon\ll 1$.
\end{enumerate}
Write $K_X+\Delta=f^*(K_C+\Lambda_C)$.
Whether there exists a positive integer $m(d, r)$ depending only on $d$ and $r$ such that $m\Lambda_C$ is base point free?
\end{question}

Note that the condition (3) above makes sure that the good moduli space $(M^K_{d,v,r}, \Lambda_{\Hodge})$ still applies to the question. The same as before, we assume all fibers admit volume $v$, then there exists a morphism $f: C\to M^K_{d,v,r}$ such that $\Lambda_C=f^*\Lambda_{\Hodge}$. The different point here is that we do not have $\fB(\Lambda_{\Hodge})\cap f(C)=\emptyset$  as in the proof of Theorem \ref{main} anymore. It is also natural to ask the same question when the base is of higher dimension.

\begin{question}\label{semiample hodge}
For any family in the K-moduli stack $\mM^K_{d,v,r}$ where the base is a projective normal variety, whether the Hodge line bundle on the base is semiample?
%Let $M$ be a component of $M^K_{d,v,r}$, whether the Hodge line bundle $\Lambda_{\Hodge}|_M$ is semiample?
\end{question}

%\begin{remark}
%This question is in fact equivalent to asking, for any component $M$ of $M^K_{d,v,r}$, whether the restriction of $\Lambda_{\Hodge}$ to $M$ is semiample.
%\end{remark}

Question \ref{semiample hodge} seems more challenging at first sight, however, it is in fact implied by Question \ref{lc fiber} (see the next subsetion), and that is why we only consider one dimensional bases.

\subsection{A criterion for semi-ampleness}

The following lemma  reduces the semi-ampleness problem to effective semi-ampleness problem on curves, which is originally due to \cite{Floris14}.

\begin{lemma}\label{reduction}
Let $X$ be a projective normal variety and $L$ a line bundle on $X$. Then the following two statements are equivalent:
\begin{enumerate}
\item $L$ is semiample,
\item there is a positive integer $m$ depending only on $(X,L)$ such that for any morphism from a projective smooth curve $f: C\to X$, the line bundle $mL_C$ is base point free. Here $mL_C$ means the pullback of $mL$.
\end{enumerate}
\end{lemma}

\begin{proof}
The direction $(1)\Rightarrow (2)$ is clear. We prove the converse direction. We use an inductive approach. Suppose for any morphism from a projective smooth variety $V\to X$ such that $\dim V\leq k$, we have that $mL_V$ is base point free.  Assume now $\dim V=k+1\geq 2$, we aim to show that $mL_V$ is base point free. 
Let $H$ be a sufficiently general very ample line bundle on $V$ such that $H-mL_V$ is ample. By Kodaira vanishing we have
$$H^1(V, mL_V-H)=0,$$
which implies the following exact sequence:
$$0\to H^0(V, mL_V-H)\to H^0(X, mL_V)\to H^0(H, mL_H)\to 0. $$
Let us write $|mL_V|=|\Mob|+\Fix$, which is the decomposition into the mobile part and fixed part. By restricting to $H$ we have the following by the above exact sequence:
$$|mL_V|_H=|mL_H|=|\Mob|_H+\Fix_H. $$
By induction, $|mL_H|$ is a free linear system, then we see that $\Fix_H=0$. This implies that the base locus of $|mL_V|$ is of codimension at least two. Let $V'\to V$ be a resolution such that the decomposition into the mobile part and fixed part
$$|mL_{V'}|=|\Mob'|+\Fix' $$
satisfies that the mobile part $|\Mob'|$ is free. For the morphism $V'\to X$ which is the composition of $V'\to V$ and $V\to X$, we have shown by induction that the linear system $|mL_{V'}|$ admits base locus of codimension at least two, thus $\Fix'=0$ and $mL_{V'}$ is base point free.  Therefore, $mL_V$ is also base point free. The proof is finished.
\end{proof}

\begin{corollary}
A positive answer to Question \ref{lc fiber} leads to a positive answer to Question \ref{semiample hodge}.
\end{corollary}

\begin{proof}
Suppose Question \ref{lc fiber} admits a positive answer. Assume we are given a family $(\mX', \frac{1}{r}\mD')\to V$ satisfying the following conditions:
\begin{enumerate}
\item the family $(\mX', \frac{c_0}{r}\mD')\to V$ is contained in $\mM^K_{d,v,r}(V)$ (see Section \ref{sec: 3} for the notation),
\item the base $V$ is a projective normal  variety.
\end{enumerate}
The aim is to show that $\Lambda_V=f^*\Lambda_{\Hodge}$  is semiample.  For any morphism $C\to V$ from a projective smooth curve $C$,  consider the family over $C$ via pullback. 
By our assumption, there is a positive integer $m$ depending only on $d$ and $r$ such that $m\Lambda_C$ is base point free. By Lemma \ref{reduction}, $\Lambda_V$ is semiample. The proof is finished.
\end{proof}

\bibliography{reference.bib}
\end{document}